\documentclass[12pt]{amsart}

\newtheorem{theorem}{Theorem}[section]
\newtheorem{lemma}[theorem]{Lemma}
\newtheorem{corollary}[theorem]{Corollary}
\theoremstyle{definition}   

\newtheorem{example}[theorem]{Example}

\theoremstyle{remark}
\newtheorem{remark}[theorem]{Remark}

\numberwithin{equation}{section}

\usepackage{times}
\usepackage{enumerate}
\usepackage{tfrupee}
\usepackage{tikz}
\usetikzlibrary{chains,fit}
\usetikzlibrary{shapes,snakes}
\usetikzlibrary{graphs}
\usepackage{graphicx,adjustbox}
\usepackage{hyperref}
\usepackage{amsmath,amssymb}

\newcommand{\lcm}{{\rm \mbox{lcm}}}
\newcommand{\LT}{{\rm \mbox{Lt}}}

\usepackage{amscd}
\usepackage{graphicx}
\usepackage[all]{xy}

\title[Ideals of the form $I_{1}(XY)$]
{Ideals of the form $I_{1}(XY)$}
\author{
Joydip Saha
\and
Indranath Sengupta
\and
Gaurab Tripathi
}
\date{}

\address{\small \rm  Discipline of Mathematics, IIT Gandhinagar, Palaj, Gandhinagar, 
Gujarat 382355, INDIA.} 
\email{saha.joydip56@gmail.com}

\address{\small \rm  Discipline of Mathematics, IIT Gandhinagar, Palaj, Gandhinagar, 
Gujarat 382355, INDIA.}
\email{indranathsg@iitgn.ac.in}
\thanks{The second author is the corresponding author.}

\address{\small \rm Department of Mathematics, Jadavpur University, Kolkata,
WB 700 032, INDIA.} 
\email{gelatinx@gmail.com}

\date{}



\begin{document}
\begin{abstract}
In this paper we compute Gr\"{o}bner bases for determinantal ideals 
of the form $I_{1}(XY)$, where $X$ and $Y$ are both matrices whose 
entries are indeterminates over a field $K$. We use the Gr\"{o}bner 
basis structure to determine Betti numbers for such ideals.
\end{abstract}

\maketitle

\section{Introduction}
Let $K$ be a field and $\{x_{ij}; \, 1\leq i \leq m, \, 1\leq j \leq n\}$, 
$\{y_{j}; \, 1\leq j \leq n\}$ be indeterminates over $K$. Let $K[x_{ij}]$ 
and $K[x_{ij}, y_{j}]$ denote the polynomial algebras over 
$K$. Let $X$ denote an $m\times n$ matrix such that its entries belong to the 
ideal $\langle \{x_{ij}; \, 1\leq i \leq m, \, 1\leq j \leq n\}\rangle$. 
Let $Y=(y_{j})_{n\times 1}$ be the generic $n\times 1$ column matrix. Let $I_{1}(XY)$ denote the ideal generated by the $1\times 1$ minors or the entries of 
the $m\times 1$ matrix $XY$. Ideals of the form $I_{1}(XY)$ appeared in 
the work of J. Herzog \cite{herzog} in 1974. These ideals are closely 
related to the notion of Buchsbaum-Eisenbud variety of complexes. A 
characteristic free study of these varieties can be found in \cite{concini}, 
where the defining equations of these varieties have been described as 
minors of matrices using combinatorial structure of multitableux. It has also 
been proved that the varieties are Cohen-Macaulay and Normal. The 
ideal $I_{1}(XY)$ is a special case of the defining ideal of a variety 
of complexes, when $n_{0} = m$, $n_{1} = n$, $n_{2} = 1$, in the notation 
of \cite{concini}. These ideals feature once again in \cite{tchernev}, 
in the study of the structure of a \textit{universal ring} of a 
\textit{universal pair} defined by Hochster. It has been proved in 
\cite{tchernev} that the set of standard monomials 
form a free basis for the universal ring. The initial ideal of the 
defining ideal is given by the set of all nonstandard monomials, which form 
a monomial ideal. A combination of Gr\"{o}bner basis techniques and 
representation theory techniques yield the results in \cite{tchernev}. We 
were not aware of this work when we computed a Gr\"{o}bner 
basis for the ideal $I_{1}(XY)$ using very elementary techniques. Our technique 
uses nothing more than the Buchberger's criterion and the description 
of Gr\"{o}bner bases for the ideals of minors of matrices from 
\cite{conca} and \cite{sturmfels}. 
\medskip

Given determinantal ideals $I$ and $J$, the sum ideal $I + J$ is often difficult to understand 
and they appear in various contexts. Ideals $I_{1}(XY)+J$ are special in the sense that they 
occur in several geometric considerations like linkage and generic residual intersection 
of polynomial ideals, especially in the context of syzygies; see \cite{nor}, \cite{akm}, 
\cite{bucheis}, \cite{bkm}, \cite{johnson}. Some important classes of ideals in this category 
are the Northcott ideals, the Herzog ideals; see Definition 3.4 in \cite{akm} and the 
deviation two Gorenstein ideals defined in \cite{hunul}. Northcott ideals were resolved by 
Northcott in \cite{nor}. Herzog gave a resolution of a special case of the Herzog ideals in 
\cite{herzog}. These results were extended in \cite{bucheis}. In a similar vein, 
Bruns-Kustin-Miller \cite{bkm} resolved the ideal $I_{1}(XY)+ I_{\min (m,n)}(X)$, where 
$X$ is a generic $m \times n$ matrix and $Y$ is a generic $n\times 1$ matrix. 
Johnson-McLoud \cite{johnson} proved certain properties for the ideals of the form 
$I_{1}(XY)+I_{2}(X)$, where $X$ is a generic symmetric matrix and $Y$ is either generic or 
generic alternating. One of the recent articles is \cite{it} which shows connection 
of ideals of this form with the ideal of the dual of the quotient bundle on the Grassmannian 
$G(2,n)$.
\medskip

Ideals of the form $I+J$ also appear naturally in the study of some natural class of curves; 
see \cite{hip}. While computing Betti numbers for such ideals, a useful technique is 
often the iterated Mapping Cone. This technique requires a good understanding 
of successive colon ideals between $I$ and $J$, which is often difficult to compute. It is helpful 
if Gr\"{o}bner bases for $I$ and $J$ are known. 
\medskip

In this paper our aim is to produce some suitable 
Gr\"{o}bner bases for ideals of the form $I_{1}(XY)$, when $Y$ is a generic column matrix and 
$X$ is one of the following:
\begin{enumerate}
\item $X$ is a generic square matrix;
\item $X$ is a generic symmetric matrix;
\item $X$ is a generic $(n+1)\times n$ matrix. 
\end{enumerate}
We have also studied $I_{1}(XY)$, when
\begin{enumerate}
\item[(4)] $X$ is an $(m\times mn)$ generic matrix and Y is an $(mn\times n)$ generic matrix.
\end{enumerate}

Our method is constructive and it would exhibit that 
the first two cases behave similarly. Newly constructed Gr\"{o}bner 
bases will be used to compute the Betti numbers of $I_{1}(XY)$. We will see that computing 
Betti numbers for $I_{1}(XY)$ in the first two cases is not difficult, while the last two 
cases are not so straightforward. We will use some results from \cite{sstprime} and 
\cite{sstsum} which have some more deep consequences of the Gr\"{o}bner basis computation 
carried out in this paper.

\section{Defining the problems}

Let $K$ be a field and $\{x_{ij}; \, 1\leq i \leq n+1, \, 1\leq j \leq n\}$, 
$\{y_{j}; \, 1\leq j \leq n\}$ be indeterminates over $K$. Let 
$R = K[x_{ij}, y_{j}\mid 1\leq i, j \leq n]$, 
$\widehat{R} = K[x_{ij}, y_{j}\mid 1\leq i \leq n+1, \, 1\leq j \leq n]$ 
denote polynomial $K$-algebras. Let $X=(x_{ij})_{n\times n}$, such that $X$ is either 
generic or generic symmetric. Let $\widehat{X}=(x_{ij})_{(n+1)\times n}$ and 
$Y=(y_{j})_{n\times 1}$ be generic matrices. We define 
$\mathcal{I} = I_{1}(XY)$ and $\mathcal{J} = I_{1}(\widehat{X}Y)$. 
\medskip

Let $g_{i} = \sum_{j=1}^{n}x_{ij}y_{j}$, for $1\leq i\leq n$. Then, $\mathcal{I} = \langle g_{1}, \ldots , g_{n}\rangle$. 
Let us choose the lexicographic monomial order on $R$ given by 
\begin{enumerate}
\item $x_{11}> x_{22}> \cdots >x_{nn}$;
\item $x_{ij}, y_{j} < x_{nn}$ for every $1 \leq i \neq j \leq n$.
\end{enumerate}
It is an interesting observation that the set $\{g_{1}, \ldots , g_{n}\}$ is 
a Gr\"obner basis for $\mathcal{I}$ with respect to the above monomial order 
and the elements $g_{1}, \ldots , g_{n}$ form a regular sequence as well; see Lemma \ref{disjoint} and Theorem \ref{regseqI}. However, this Gr\"{o}bner 
basis is too small in size to be of much help in applications like 
computing primary decomposition of $I_{1}(XY)$ or computing Betti numbers of 
ideals of the form $I_{1}(XY) + J$, carried out in \cite{sstprime} 
and \cite{sstsum} respectively. This motivated us to look for a                                        a different Gr\"{o}bner basis for $\mathcal{I}$; see Theorem \ref{gbtheorem}. 
This construction gives rise to a 
bigger picture and naturally generalizes to a Gr\"{o}bner basis for the ideal 
$\mathcal{J} = I_{1}(\widehat{X}Y)$. As an application, we compute the Betti 
numbers for the ideals $\mathcal{I}$ and $\mathcal{J}$; see section 6. 

\section{Notation}

\begin{enumerate}
\item[(i)]  $C_{k}:= \{\mathbf{a} = (a_{1},\cdots,a_{k}) \mid 1\leq  a_{1} < \cdots < a_{k} \leq n\}$; denotes the collection of all ordered $k$-tuples from $\{1,\cdots,n\}$. In case of $\mathcal{J} = I_{1}(\widehat{X}Y)$, the set $C_{k}$ would denote the collection of all ordered $k$-tuples $(a_{1},\cdots,a_{k})$ from $\{1,\cdots,n+1\}$.
\item[(ii)]  Given $\mathbf{a}=(a_{1},\ldots,a_{k})\in C_{k}$;
\begin{itemize}
\item $X^{\mathbf{a}}=[a_{1},\cdots,a_{k}|1,2,\ldots , k]$ denotes the $k\times k$ minor of the matrix $X$, with $a_{1},\ldots,a_{k}$ as rows and $1,\ldots , k$ as columns. Similarly, $\widehat{X}^{\mathbf{a}}=[a_{1},\cdots,a_{k}|1,\ldots , k]$ denotes the $k\times k$ minor of the matrix $\widehat{X}$, with 
$a_{1},\ldots,a_{k}$ as rows and $1,\ldots , k$ as columns.
\item $S_{k}:=\{X^{\mathbf{a}}:\mathbf{a}\in C_{k}\}$ and  $I_{k}$ denotes the ideal generated by $S_{k}$ in the polynomial ring $R$ (respectively $\widehat{R}$);
\item $X^{\mathbf{a},m}:= [a_{1},\cdots,a_{k}|1,\cdots,k-1,m]$ if $m\geq k$;
\item $\widetilde{X^{\mathbf{a}}} 
= \sum_{m\geq k}[a_{1},\cdots,a_{k}|1,\cdots,k-1,m]y_{m} 
= \sum_{m\geq k}X^{\mathbf{a},m}y_{m}$;
\item $\widetilde{S}_{k}:= \{\widetilde{X^{\mathbf{a}}}: X^{\mathbf{a}} \in S_{k}\}$ and $\widetilde{I}_{k}$ denotes the ideal generated by 
$\widetilde{S}_{k}$ in the polynomial ring $R$ (respectively $\widehat{R}$);
\item  $G_{k}=\cup_{i\geq k}\widetilde{S}_{i}$;
\item  $G =\cup_{k\geq 1}G_{k}$;
\item $X^{\mathbf{a}}_{r}:= 
[a_{1},a_{2},\cdots,\hat{a_{r}},a_{r+1}\cdots,a_{k}|1,2,\cdots,k-1]$, 
if $k\geq 2$.
\end{itemize} 
\item[(iii)] Suppose that $C_{k} = \left\lbrace\mathbf{a}_{1}< \ldots < \mathbf{a}_{\binom{n}{k}}\right\rbrace$, where $<$ is the lexicographic 
ordering. Given $m\geq k$, the map 
$$\sigma_{m}:\left\lbrace X^{\mathbf{a}_{1},m}, \ldots , X^{\mathbf{a}_{\binom{n}{k}},m}\right\rbrace \rightarrow 
\left\lbrace 1,\cdots,\binom{n}{k}\right\rbrace$$ is defined by $\sigma_{m}(X^{\mathbf{a}_{i},m})=i$. This is a bijective map. The map $\sigma_{k}$ will be denoted by $\sigma$, which is the bijection from $S_{k}$ 
to $\{1,\cdots,\binom{n}{k}\}$ given by $\sigma(X^{\mathbf{a}_{i}})=\sigma_{k}(X^{\mathbf{a}_{i}, k})=i$.
\end{enumerate}

\section{Gr\"{o}bner basis for $\mathcal{I}$}
We first construct a Gr\"{o}bner basis for the ideal $\mathcal{I}$. A similar computation 
works for computing a Gr\"{o}bner basis for the ideal $\mathcal{J}$, which will be discussed in the next section. Our aim in this section is to prove 

\begin{theorem}\label{gbtheorem}
The set $G_{k}$ is a reduced Gr\"{o}bner Basis for the ideal $\widetilde{I}_{k}$, 
with respect to the lexicographic monomial order induced by the following order 
on the variables: $y_{1}>y_{2}>\cdots >y_{n}> x_{ij}$ for all $i,j$, such 
that $x_{ij}>x_{i'j'}$ if $i<i'$ or if $i=i'$ and $j< j'$. In 
particular, $\mathcal{G} = G_{1}$ is a reduced Gr\"{o}bner Basis for the ideal 
$\widetilde{I}_{1} = \mathcal{I}$.
\end{theorem}
\medskip

\noindent We first write down the main steps involved in the proof. Let 
$\widetilde{X}^{\mathbf{a}},\widetilde{X}^{\mathbf{b}}\in G_{k}=\cup_{i\geq k}\widetilde{S}_{i}$. Then, either 
$X^{\mathbf{a}},X^{\mathbf{b}}\in S_{k}$ or $X^{\mathbf{a}}\in S_{k}$, $X^{\mathbf{b}}\in S_{k'}$, for $k'>k$. Our aim is to show that 
$S(\widetilde{X}^{\mathbf{a}},\widetilde{X}^{\mathbf{b}})\rightarrow_{G_{k}} 0$ and use Buchberger's criterion.
\begin{enumerate}
\item[(A)] By Lemma \ref{concath}, we have $S(X^{\mathbf{a}},X^{\mathbf{b}})\longrightarrow_{S_{k}} 0$. We write 
$m_{\mathbf{a}}X^{\mathbf{a}}+m_{\mathbf{b}}X^{\mathbf{b}}
=S(X^{\mathbf{a}},X^{\mathbf{b}})
= \sum_{t=1}^{\binom{n}{k}}\alpha_{t}X^{\mathbf{a}_{t}}\longrightarrow_{S_{k}} 0$, such that 
$X^{\mathbf{a}_{i}} = X^{\mathbf{a}}$ and $X^{\mathbf{a}_{j}} = X^{\mathbf{b}}$, for some $i$ and $j$. Therefore, by Schreyer's theorem the tuples 
$(\alpha_{1},\ldots ,\alpha_{i}-m_{\mathbf{a}},\ldots , \alpha_{j}-m_{\mathbf{b}},\ldots,\alpha_{r} )$ generate 
${\rm Syz}(I_{k})$. 
\item[(B)] ${\rm Syz}(I_{k})$ is precisely known by \cite{en}.
\item[(C)] $S(\widetilde{X}^{\mathbf{a}},\widetilde{X}^{\mathbf{b}})\longrightarrow_{\widetilde{S}_{k}} S(\widetilde{X}^{\mathbf{a}},\widetilde{X}^{\mathbf{b}}) 
- \sum_{t=1}^{\binom{n}{k}} \alpha_{t}\widetilde{X}^{\mathbf{a}_{t}}$ by Lemma \ref{spolylemma}, if $X^{\mathbf{a}},X^{\mathbf{b}}\in S_{k}$ and by 
Lemma \ref{spolydifferent}, if $X^{\mathbf{a}}\in S_{k}$, $X^{\mathbf{b}}\in S_{k'}$, for $k'>k$.
\item[(D)] $S(\widetilde{X}^{\mathbf{a}},\widetilde{X}^{\mathbf{b}})- \sum_{t=1}^{\binom{n}{k}} \alpha_{t}\widetilde{X}^{\mathbf{a}_{t}}=s\in \widetilde{I}_{k+1}$, by Lemma \ref{spolylemma}, if $X^{\mathbf{a}},X^{\mathbf{b}}\in S_{k}$.
\item[(E)] $S(\widetilde{X}^{\mathbf{a}},\widetilde{X}^{\mathbf{b}})- \sum_{t=1}^{\binom{n}{k}} \alpha_{t}\widetilde{X}^{\mathbf{a}_{t}}=s\in \widetilde{I}_{k'+1}$, by Lemma \ref{spolydifferent}, if $ X^{\mathbf{a}}\in S_{k}$, $X^{\mathbf{b}}\in S_{k'}$, for $k'>k$.
\item[(F)] $s\longrightarrow_{G_{k}} 0$, proved in Theorem \ref{gbtheorem} for both the cases.
\end{enumerate}
\medskip

We first prove a number of Lemmas to complete the proof through the steps mentioned 
above. 
\medskip

\begin{lemma}\label{concath}
The set $S_{k}$ forms a Gr\"{o}bner basis of $I_{k}$ with respect to the chosen 
monomial order on $R$.
\end{lemma}

\proof We use Buchberger's criterion for the proof. Let $\mathbf{c}, \mathbf{d}\in S_{k}$. Suppose that $S(X^{\mathbf{c}},X^{\mathbf{d}})\stackrel{S_{k}}{\longrightarrow} r$. Then,  $S(X^{\mathbf{c}},X^{\mathbf{d}})-\sum_{\mathbf{a_{i}}\in C_{i}}h_{i}X^{\mathbf{a_{i}}}=r$. 
\medskip
 
If $X$ is generic (respectively generic symmetric), we know by \cite{sturmfels}
(respectively by \cite{conca}) that the set of all $k\times k$ minors of the 
matrix $X$ forms a Gr\"{o}bner basis for the ideal $I_{k}(X)$, with respect 
to the chosen monomial order. Therefore, there exists $[a_{1},a_{2},\cdots,a_{k}\mid b_{1},b_{2},\cdots,b_{k}]$, such that its leading term $\prod_{i=1}^{k}x_{a_{i}b_{i}}$ divides $\LT(r)$. We see that if $b_{k}= k$, the minor belongs 
to the set $S_{k}$ and we are done. 
\medskip

Let us now consider the case $b_{k}\geq k+1$. Let $X$ be generic symmetric. Then, $a_{k} = k$ and $b_{k}\geq k+1$ imply that the minor belongs to the set $S_{k}$. If $a_{k},b_{k}\geq k+1$, then $x_{a_{k}b_{k}}\mid \LT(r)$ but $x_{a_{k}b_{k}}$ doesn't divide any term of elements in $S_{k}$. Let $X$ be generic. Then, for 
any $a_{k}$ and under the condition $b_{k}\geq k+1$, then $x_{a_{k}b_{k}}\mid \LT(r)$ but $x_{a_{k}b_{k}}$ doesn't divide any term of elements in $S_{k}$.\qed
\medskip

\begin{lemma}\label{disjoint}
Let $h_{1},h_{2}\cdots, h_{n}\in R$ be such that with respect to a suitable 
monomial order on $R$, the leading terms of them are pairwise coprime. 
Then, $h_{1},h_{2}\cdots, h_{n}$ is a Gr\"{o}bner basis of the ideal 
generated by $h_{1},h_{2}\cdots, h_{n}$ with respect to the same monomial 
order and they form a regular sequence in $R$.
\end{lemma}

\proof. The proof is a routine application of the division algorithm.\qed
\medskip

\begin{lemma}\label{EN}
Let $1\leq k \leq n$. The height of the ideal $I_{k}$ is $n-k+1$, in case of $X$. 
\end{lemma}

\proof. Let us consider the case for $X$. We know that $ht(I_{k})\leq n-k+1$. 
It suffices to find a regular sequence of that length in the ideal $I_{k}$. 
We claim that $\{[1\cdots k|1 \cdots k], [2\cdots k+1|1\cdots k],\ldots ,
[n-k+1 \cdots n|1 \cdots k]\}$ forms a regular sequence. The leading term of $[a_{1},a_{2},\cdots,a_{k}\mid b_{1},b_{2},\cdots,b_{k}]$ with respect to the chosen monomial order is $\prod_{i=1}^{k}x_{a_{i}b_{i}}$. Therefore, leading terms of the above minors are mutually coprime and we are done by Lemma \ref{disjoint}. \qed
\medskip

\begin{remark}
We now assume that $X = (x_{ij})$ is a generic $n\times n$ matrix. The proof 
for the symmetric case is exactly the same.
\end{remark}
\medskip
 
\subsection*{Description of generators of ${\rm Syz}(I_{k})$} 
By Lemma \ref{EN} we conclude that a minimal free resolution of the ideal $I_{k}$ 
is given by the Eagon-Northcott complex. Let us describe the first syzygies of 
the Eagon-Northcott resolution of $I_{k}$.
\medskip

Let $\mathbf{a} = (a_{1}, \ldots , a_{k+1})\in C_{k+1}$. For $1\leq r\leq k+1$, we define $X^{\mathbf{a}}_{r}=[a_{1},\ldots,\hat{a_{r}}, \ldots,a_{k+1}|1,\ldots,k]$.  Hence 
$X^{\mathbf{a}}_{r}\in S_{k}$. We define the map $\phi$ as follows.
\begin{eqnarray*}
\{1,2,\cdots,k\}\times C_{k+1}& \stackrel{\phi}{\longrightarrow} & R^{\binom{n}{k}}\\
(j,\mathbf{a})& \mapsto & \alpha
\end{eqnarray*}
\medskip

such that $\, \alpha(i) = 
\begin{cases}
(-1)^{r_{i}+1}x_{(a_{r_{i}},\ j)} & \textrm{if} \, i=\sigma(X^{\mathbf{a}}_{r_{i}})\ \textrm{for \ some} \quad r_{i}; \\[2mm]
0 & \textrm{otherwise}.
\end{cases}$
\medskip

\noindent The map $\sigma$ is the bijection from $S_{k}$ to $\{1,2,\cdots,\binom{n}{k}\}$, defined 
before. The image of $\phi$ gives a complete list of generators of ${\rm Syz}(I_{k})$.
\medskip

\begin{example}
We give an example, by taking $k=3$ and $n=5$. Let 
$\sigma: S_{5} \longrightarrow \{1,\cdots \binom{5}{3}\}$ be defined by,
\begin{itemize}
\item $[1,2,3\mid 1,2,3]\mapsto 1$
\item $[1, 2, 4\mid 1,2,3]\mapsto 2$
\item $[1,2,5\mid 1,2,3]\mapsto 3$
\item $[1, 3, 4\mid 1,2,3]\mapsto 4$
\item $[1, 3, 5\mid 1,2,3]\mapsto 5$
\item $[1,4,5\mid 1,2,3]\mapsto 6$
\item $[2,3,4\mid 1,2,3]\mapsto 7$
\item $[2,3,5\mid 1,2,3]\mapsto 8$
\item $[2,4,5\mid 1,2,3]\mapsto 9$
\item $[3,4,5\mid 1,2,3]\mapsto 10$
\end{itemize} 
In our example, $\phi:\{1,\cdots 3\}\times C_{4}\longrightarrow R^{\binom{5}{3}}$ and 
$\phi (j,\mathbf{a})\mapsto \alpha$. Let $j=2$ and $\mathbf{a}=(1, 3, 4, 5)$. Then, $X^{\mathbf{a}}_{1}=[3, 4, 5\mid 1, 2, 3]$, $X^{\mathbf{a}}_{2}=[1, 4, 5\mid 1, 2, 3]$, 
$X^{\mathbf{a}}_{3}=[1, 3, 5\mid 1, 2, 3]$, $X^{\mathbf{a}}_{4}=[1, 3, 4\mid 1, 2, 3]$. 
Therefore, $\sigma(X^{\mathbf{a}}_{1})=10$, 
$\sigma(X^{\mathbf{a}}_{2})=6$, 
$\sigma(X^{\mathbf{a}}_{3})=5$, 
$\sigma(X^{\mathbf{a}}_{4})=4$. Similarly, 
$\alpha(4)=(-1)^{4+1}x_{52}=-x_{52}$, $\alpha(5)=(-1)^{3+1}x_{42}=x_{42}$, 
$\alpha(6)=(-1)^{2+1}x_{32}=-x_{32}$, $\alpha(10)=(-1)^{1+1}x_{12}=x_{12}$. 
Therefore, $\alpha=(0,0,0,-x_{52},x_{42},-x_{32},0,0,0,x_{12})$.\\
\end{example}

\begin{lemma}\label{syzlemma} 
Let $1\leq k\leq n-1$ and let $S_{k}=\left\lbrace X^{\mathbf{a}_{1}}, \ldots , X^{\mathbf{a}_{\binom{n}{k}}}\right\rbrace$ be such that $\mathbf{a}_{1}< \ldots < \mathbf{a}_{\binom{n}{k}}$ with respect to the lexicographic ordering. Suppose that $\alpha =(\alpha_{1},\cdots,\alpha_{\binom{n}{k}})\in$ Syz$^{1}(I_{k})$, then $\sum_{i=1}^{\binom{n}{k}}\alpha_{i}X^{\mathbf{a}_{i}}=0$ and   $\sum_{i=1}^{\binom{n}{k}}\alpha_{i}\widetilde{X^{\mathbf{a}_{i}}}\in \widetilde{I}_{k+1}$. 
\end{lemma}

\proof  
We have $\widetilde{X}^{\mathbf{a}_{i}}=\sum_{m\geq k}\sigma_{m}^{-1}(i)y_{m}$. Therefore 
$$\sum_{i=1}^{\binom{n}{k}}\alpha_{i}\widetilde{X}^{\mathbf{a}_{i}}=\sum_{i}\alpha_{i}(\sum_{m\geq k}\sigma_{m}^{-1}(i)y_{m})=\sum_{m\geq k}(\sum_{i}\alpha_{i}\sigma_{m}^{-1}(i)y_{m}).$$ 
It is enough to show that $\sum_{i}\alpha_{i}\sigma_{m}^{-1}(i)y_{m}\in \widetilde{I}_{k+1}$, 
for every $m\geq k$. We have $\alpha\in {\rm Syz}(I_{k})=\langle {\rm Im}(\phi) \rangle$. 
Without loss of generality we may assume that 
$\alpha\in {\rm Im}(\phi)$. There exists $(j,\mathbf{a}_{k+1})\in \{1,2,\cdots k\}\times C_{k+1}$ such that $\phi(j,\mathbf{a}_{k+1})=\alpha$. We will show that 
$\alpha_{i}\cdot \sigma_{m}^{-1}(i)\in I_{k+1}$ for every $m\geq k$ and each $i$. We have 
$i=\sigma (X_{r_{i}}^{\mathbf{a}_{k+1}})$ since $\alpha_{i}\neq 0$. 
But $\sigma_{m}^{-1}(i)=[a_{1},\ldots,\hat{a_{r_{i}}},\ldots,a_{k+1}|1,\ldots,k-1,m ]$. 
We have 
$$[a_{1},\ldots,a_{k+1}|j,1,\ldots,k-1,m]=0\quad \text{for} \quad j\leq k-1\quad \text{and}$$ 
$$[a_{1},\ldots,a_{k+1}|k,1,\ldots,k-1,m]= (-1)^{k}[a_{1},\ldots,a_{k+1}|1,\ldots,k,m]\in I_{k+1}.$$ 
Therefore, 
\begin{eqnarray*}
\sum_{i=1}^{\binom{n}{k}}\alpha_{i}\cdot\sigma_{m}^{-1}(i) & = & \sum_{i=1}^{\binom{n}{k}} (-1)^{r_{i}+1}x_{(a_{r_{i}},\ j)}[a_{1},\ldots,\hat{a_{r_{i}}},\ldots,a_{k+1}|1,\ldots,k-1,m]\\
{} & = & [a_{1},\ldots,a_{k+1}|j,1,\ldots,k-1,m]\in I_{k+1};
\end{eqnarray*}
Hence, 
$$\sum_{i=1}^{\binom{n}{k}}\alpha_{i}\widetilde{X^{\mathbf{a}_{i}}} =\sum_{i=1}^{\binom{n}{k}}\alpha_{i}\cdot \widetilde{\sigma_{m}^{-1}(i)}= 
(-1)^{k}\sum_{i=1}^{\binom{n}{k}}[a_{1},\ldots,a_{k+1}|1,\ldots,k,m]y_{m}\in \widetilde{I}_{k+1}.\qed$$

\begin{lemma}\label{spolylemma}
Let $X^{\mathbf{a}_{i}},X^{\mathbf{a}_{j}}\in S_{k}=\left\lbrace X^{\mathbf{a}_{1}}, \ldots , X^{\mathbf{a}_{\binom{n}{k}}}\right\rbrace$, for $i\neq j$. Then, there exist monomials $h_{t}$ in $R$ and a polynomial $r\in \widetilde{I}_{k+1}$ such that 
\begin{enumerate}
\item[(i)] $S(X^{\mathbf{a}_{i}},X^{\mathbf{a}_{j}})
=\sum_{t=1}^{\binom{n}{k}} h_{t}X^{\mathbf{a}_{t}}$, upon division 
by $S_{k}$;
\item[(ii)] $S(\widetilde{X}^{\mathbf{a}_{i}},\widetilde{X}^{\mathbf{a}_{j}})=\sum_{t=1}^{\binom{n}{k}} h_{t}\widetilde{X}^{\mathbf{a}_{t}}+r$, upon division by $\widetilde{S}_{k}$. 
\end{enumerate}
\end{lemma} 

\proof (i) The expression follows from the observation that $S_{k}$ is a Gr\"{o}bner 
basis for the ideal $I_{k}$.
\medskip

\noindent (ii) We first note that, $\LT(\widetilde{X}^{\mathbf{a}_{t}})=\LT(X^{\mathbf{a}_{t}})y_{k}$, for 
every $X^{\mathbf{a}_{t}}\in S_{k}$.
Let $S(X^{\mathbf{a}_{i}},X^{\mathbf{a}_{j}})=cX^{\mathbf{a}_{i}}-dX^{\mathbf{a}_{j}}$, 
where $c=\dfrac{\lcm(\LT(X^{\mathbf{a}_{i}}),\LT(X^{\mathbf{a}_{j}}))}{X^{\mathbf{a}_{i}}}$ and 
$d=\dfrac{\lcm(\LT(X^{\mathbf{a}_{i}}),\LT(X^{\mathbf{a}_{j}}))}{X^{\mathbf{a}_{j}}}$

Hence, 
\begin{eqnarray*}
S(\widetilde{X}^{\mathbf{a}_{i}},\widetilde{X}^{\mathbf{a}_{j}}) & = &
c\cdot \widetilde{X}^{\mathbf{a}_{i}}-d\cdot \widetilde{X}^{\mathbf{a}_{i}}\\
{} & = & \sum_{m\geq k}\left[c\cdot X^{\mathbf{a}_{i}, m}-d\cdot X^{\mathbf{a}_{j}, m}\right]y_{m}.
\end{eqnarray*}
It follows immediately that $\LT(S(\widetilde{X}^{\mathbf{a}_{i}},\widetilde{X}^{\mathbf{a}_{j}}))=y_{k}\LT(S(X^{\mathbf{a}_{i}},X^{\mathbf{a}_{j}}))$.
\medskip

The set $S_{k}$ is a Gr\"{o}bner basis for the ideal $I_{k}$. Therefore, 
we have $\LT(X^{\mathbf{a}_{t}})\mid \LT(S(X^{\mathbf{a}_{i}},X^{\mathbf{a}_{j}}))$, 
for some $t$. Then, $\LT(\widetilde{X}^{\mathbf{a}_{t}})\mid \LT(S(\widetilde{X}^{\mathbf{a}_{i}},\widetilde{X}^{\mathbf{a}_{j}}))$ and we have $h_{t}=\dfrac{\LT(S(X^{\mathbf{a}_{i}},X^{\mathbf{a}_{j}}))}{\LT(X^{\mathbf{a}_{t}})}= \dfrac{\LT(S(\widetilde{X}^{\mathbf{a}_{i}},\widetilde{X}^{\mathbf{a}_{j}}))}{\LT(\widetilde{X}^{\mathbf{a}_{t}})}$. We can write
\begin{eqnarray*}
r_{1} & := & S(\widetilde{X}^{\mathbf{a}_{i}},\widetilde{X}^{\mathbf{a}_{j}})- h_{t}\widetilde{X}^{\mathbf{a}_{t}}\\
{} &  = & \sum_{m\geq k}[c\cdot X^{\mathbf{a}_{i},m}-d\cdot X^{\mathbf{a}_{j},m}-h_{t}X^{\mathbf{a}_{t},m}]y_{m}\\
{} & = & \sum_{m> k}[c\cdot X^{\mathbf{a}_{i},m}-d\cdot X^{\mathbf{a}_{j},m}-h_{t}X^{\mathbf{a}_{t},m}]y_{m}+[c\cdot X^{\mathbf{a}_{i}}-d\cdot X^{\mathbf{a}_{j}}-h_{t}X^{\mathbf{a}_{t}}] y_{k}
\end{eqnarray*}
Note that $r_{1}\in \widetilde{I}_{k}$ and 
$\LT(r_{1}) = \LT (S(\widetilde{X}^{\mathbf{a}_{i}},\widetilde{X}^{\mathbf{a}_{j}})- h_{t}\widetilde{X}^{\mathbf{a}_{t}})= y_{k}\LT(S(X^{\mathbf{a}_{i}},X^{\mathbf{a}_{j}})-h_{t}X^{\mathbf{a}_{t}})$. We proceed as before with the polynomial 
$S(X^{\mathbf{a}_{i}},X^{\mathbf{a}_{j}})-h_{t}X^{\mathbf{a}_{t}}\in I_{k}$ and continue the 
process to obtain the desired expression involving the polynomial $r$.
\medskip

We now show that the polynomial $r$ is in the ideal $\widetilde{I}_{k+1}$. 
Let us write $H_{j}= h_{j}+ d$, $H_{i}= h_{i}-c$ 
and $H_{t}= h_{t}$ for $t\neq i,j$. It follows from $S(X^{\mathbf{a}_{i}},X^{\mathbf{a}_{j}})=
\sum_{t=1}^{\binom{n}{k}} h_{t}X^{\mathbf{a}_{t}}$, that $\sum_{t=1}^{\binom{n}{k}}H_{t}X^{\mathbf{a}_{t}}= 0$. Therefore, 
$\mathbf{H} = (H_{1}, \ldots , H_{\binom{n}{k}})\in {\rm Syz}(I_{k})$ and 
by Lemma \ref{syzlemma} we have 
$\sum_{t=1}^{\binom{n}{k}}H_{t}\widetilde{X}^{\mathbf{a}_{t}}\in \widetilde{I}_{k+1}$. 
Hence, $r = S(\widetilde{X}^{\mathbf{a}_{i}},\widetilde{X}^{\mathbf{a}_{j}})-\sum_{t\neq i,j} h_{t}\widetilde{X}^{\mathbf{a}_{t}}\in \widetilde{I}_{k+1}$. \qed
\medskip

\begin{lemma}\label{laplace}
\begin{enumerate}  
\item [(i)] Let $k^{'}>k$ and $\mathbf{a} = (a_{1}, \ldots , a_{k^{'}})\in C_{k^{'}}$. 
Suppose that 
$X^{\mathbf{a}}=\sum_{\mathbf{b}_{t}\in C_{k}}\beta_{\mathbf{b}_{t}}X^{\mathbf{b}_{t}}$ is the Laplace expansion of $X^{\mathbf{a}}$. 
Then 
$$\sum_{\mathbf{b}_{t}\in C_{k}}\beta_{\mathbf{b}_{t}}X^{\mathbf{b}_{t},i} = [a_{1}, \ldots , a_{k^{'}}|1, \ldots , k-1, i, k+1, \ldots , k^{'}].$$
\item [(ii)] Let $k^{'}>k$; $\mathbf{a}=(a_{1},\ldots,a_{k^{'}})\in C_{k^{'}}$, 
$\mathbf{b}=(b_{1},\ldots,b_{k})\in C_{k}$. Suppose that 
$X^{\mathbf{a}}=\sum_{\mathbf{p}\in C_{k}}\alpha_{\mathbf{p}}X^{\mathbf{p}}$ 
and $S(X^{\mathbf{a}},X^{\mathbf{b}})=cX^{\mathbf{a}}-dX^{\mathbf{b}}=\sum_{\mathbf{p}\in C_{k}}\beta_{\mathbf{p}}X^{\mathbf{p}}$. Then 
$$\quad c\sum_{t\geq k}[a_{1}, \cdots , a_{k^{'}}| 1, \cdots, k-1, t, k+1, \cdots , k^{'}]y_{t}-d\widetilde{X}^{\mathbf{b}}-\sum_{\mathbf{p}\in C_{k}}\beta_{\mathbf{p}}\widetilde{X}^{\mathbf{p}}\in \widetilde{I}_{k+1}.$$
\end{enumerate}
\end{lemma}

\proof (i) See \cite{janjic}.
\medskip

\noindent (ii) We have $S(X^{\mathbf{a}},X^{\mathbf{b}})=cX^{\mathbf{a}}-dX^{\mathbf{b}}=\sum_{\mathbf{p}\in C_{k}}\beta_{\mathbf{p}}X^{\mathbf{p}}$. By rearranging terms we get 
$\sum_{\mathbf{p}\in C_{k}}(c\alpha_{\mathbf{p}}-\beta_{\mathbf{p}})X^{\mathbf{p}}-dX^{\mathbf{b}}=0$ and by separating out the term $(c\alpha_{\mathbf{b}}-\beta_{\mathbf{b}})X^{\mathbf{b}}$ we get 
$\sum_{\mathbf{p}\neq\mathbf{b}}(c\alpha_{\mathbf{p}}-\beta_{\mathbf{p}})X^{\mathbf{p}}+(c\alpha_{\mathbf{b}}-\beta_{\mathbf{b}}-d)X^{\mathbf{b}}=0$. Therefore, 
$\sum_{\mathbf{p}\neq\mathbf{b}}(c\alpha_{\mathbf{p}}-\beta_{\mathbf{p}})\widetilde{X}^{\mathbf{p}}+(c\alpha_{\mathbf{b}}-\beta_{\mathbf{b}}-d)\widetilde{X}^{\mathbf{b}}
\in \widetilde{I}_{k+1}$, by Lemma \ref{syzlemma}. Hence $\sum_{t\geq k}\sum_{\mathbf{p}\neq\mathbf{b}}(c\alpha_{\mathbf{p}}-\beta_{\mathbf{p}})X^{\mathbf{p},t}y_{t}+(c\alpha_{\mathbf{b}}-\beta_{\mathbf{b}}-d)\sum_{t\geq k}X^{\mathbf{b},t}y_{t}\in \widetilde{I}_{k+1}$. Now $\quad \sum_{t\geq k}\sum_{\mathbf{p}\in C_{k}}\alpha_{\mathbf{p}}X^{p,t}=  \quad \sum_{t\geq k}[a_{1}, \cdots , a_{k^{'}}| 1, \cdots, k-1, t, k+1, \cdots , k^{'}]$ by (i). Hence, 
$$c\sum_{t\geq k}[a_{1}, \cdots , a_{k^{'}}| 1, \cdots, k-1, t, k+1, \cdots , k^{'}]y_{t}-d\widetilde{X}^{\mathbf{b}}-\sum_{\mathbf{p}\in C_{k}}\beta_{\mathbf{p}}\widetilde{X}^{\mathbf{p}}\in \widetilde{I}_{k+1}.\qed $$ 

\begin{lemma}\label{spolydifferent}
Let $k^{'}>k$; $\mathbf{a}=(a_{1},\ldots,a_{k^{'}})\in C_{k^{'}}$, 
$\mathbf{b}=(b_{1},\ldots,b_{k})\in C_{k}$. Suppose that 
$S_{k}=\left\lbrace X^{\mathbf{a}_{1}}, \ldots , X^{\mathbf{a}_{\binom{n}{k}}}\right\rbrace $, such that $\mathbf{a}_{1}< \ldots < \mathbf{a}_{\binom{n}{k}}$ with respect to the lexicographic ordering. Then, there exist monomials 
$h_{t}\in R$ and a polynomial $r\in \widetilde{I}_{k+1}$ such that 
\begin{enumerate}
\item [(i)] $S(X^{\mathbf{a}},X^{\mathbf{b}})=\sum_{t=1}^{\binom{n}{k}}h_{t}X^{\mathbf{a}_{t}}$, upon division by $S_{k}$.
\item [(ii)] $S(\widetilde{X}^{\mathbf{a}},\widetilde{X}^{\mathbf{b}})=\sum_{t=1}^{\binom{n}{k}}(h_{t}\widetilde{X}^{\mathbf{a}_{t}})y_{k'} + r$, upon division by $\widetilde{S}_{k}$.
\end{enumerate}
\end{lemma}
\proof (i) The expression follows from the observation that $S_{k}$ is a Gr\"{o}bner 
basis for the ideal $I_{k}$. 
\medskip

\noindent (ii) Let $S(X^{\mathbf{a}},X^{\mathbf{b}})=cX^{\mathbf{a}}-dX^{\mathbf{b}}$, 
where $c=\dfrac{\lcm (\LT(X^{\mathbf{a}}),\LT(X^{\mathbf{b}}))}{X^{\mathbf{a}}}$ and 
$d=\dfrac{\lcm (\LT(X^{\mathbf{a}}),\LT(X^{\mathbf{b}}))}{X^{\mathbf{b}}}$. Then,
\begin{eqnarray*}
S(\widetilde{X}^{\mathbf{a}},\widetilde{X}^{\mathbf{b}}) & = & cy_{k}\widetilde{X}^{\mathbf{a}}-dy_{k^{'}}\widetilde{X}^{\mathbf{b}}\\
{} & = & 
cy_{k}\sum_{t\geq k^{'}}X^{\mathbf{a},t}y_{t}-dy_{k^{'}}\sum_{t\geq k}X^{\mathbf{b},t}y_{t}\\
{} & = & y_{k}y_{k^{'}}(cX^{\mathbf{a}}-dX^{\mathbf{b}})+ \text{terms \,devoid \,of} \, y_{k}.
\end{eqnarray*}
We therefore have 
$\text{Lt}(S(\widetilde{X}^{\mathbf{a}},\widetilde{X}^{\mathbf{b}}))=y_{k}y_{k^{'}}\text{Lt}(S(X^{\mathbf{a}},X^{\mathbf{b}}))$, since $y_{k}$ is the largest variable appearing 
in the above expression. The set $S_{k}$ being a Gr\"{o}bner basis for the ideal $I_{k}$, we have $\LT(X^{\mathbf{a_{t}}})$ dividing 
$\LT(S(X^{\mathbf{a_{i}}},X^{\mathbf{a_{j}}}))$ for some $t$. Let 
$h_{t}=\dfrac{\LT(cX^{\mathbf{a}}-dX^{\mathbf{b}})}{\LT(X^{\mathbf{a}_{t}})}$, 
with $t=1, \ldots , \binom{n}{k}$. Moreover, $\LT(\widetilde{X}^{\mathbf{a}_{t}})$ being equal to 
$y_{k}\LT(X^{\mathbf{a}_{t}})$, it divides 
$\LT(S(\widetilde{X}^{\mathbf{a}},\widetilde{X}^{\mathbf{b}}))$. Let  
$$r_{1}:= S(\widetilde{X}^{\mathbf{a}},\widetilde{X}^{\mathbf{b}})-\frac{\LT(S(\widetilde{X}^{\mathbf{a}},\widetilde{X}^{\mathbf{b}}))}{\LT(\widetilde{X}^{\mathbf{a}_{t}})}\widetilde{X}^{\mathbf{a}_{t}}=S(\widetilde{X}^{\mathbf{a}},\widetilde{X}^{\mathbf{b}})-y_{k^{'}}h_{t}\widetilde{X}^{\mathbf{a}_{t}}\in\widetilde{I}_{k}.$$
We have
\begin{eqnarray*}
r_{1} & = & y_{k}y_{k^{'}}(cX^{\mathbf{a}}-dX^{\mathbf{b}})-y_{k^{'}}h_{t}\widetilde{X}^{\mathbf{a}_{t}}+ \text{terms \,devoid \,of} \, y_{k}\\
{} & = & y_{k}y_{k^{'}}(cX^{\mathbf{a}}-dX^{\mathbf{b}}) - y_{k^{'}}h_{t}\sum_{i\geq k}X^{\mathbf{a}_{t},i}y_{i}+ \text{terms \,devoid \,of} \, y_{k}\\
{} & = & y_{k}y_{k^{'}}(cX^{\mathbf{a}}-dX^{\mathbf{b}} - h_{t}X^{\mathbf{a}_{t}})+ \text{terms \,devoid \,of} \, y_{k}\\
{} & = & y_{k}y_{k^{'}}(S(X^{\mathbf{a}},X^{\mathbf{b}}) - h_{t}X^{\mathbf{a}_{t}}) + \text{terms \,devoid \,of} \, y_{k}.\\
\end{eqnarray*}
Hence, $\LT(r_{1}) = \LT (S(X^{\mathbf{a}},X^{\mathbf{b}}) - h_{t}X^{\mathbf{a}_{t}})
= y_{k}y_{k'}\LT(S(X^{\mathbf{a}},X^{\mathbf{b}}) - h_{t}X^{\mathbf{a}_{t}})$. We proceed 
as before with the polynomial 
$S(X^{\mathbf{a}},X^{\mathbf{b}})-h_{t}X^{\mathbf{a}_{t}}\in I_{k}$ 
and continue the process to obtain the desired expression involving 
the polynomial $r$.
\medskip

We now show that the polynomial $r$ is in the ideal $\widetilde{I}_{k+1}$. 
Let us write  
\begin{eqnarray*}
r & = & S(\widetilde{X}^{\mathbf{a}},\widetilde{X}^{\mathbf{b}})-\sum_{t=1}^{\binom{n}{k}}(h_{t}\widetilde{X}^{\mathbf{a}_{t}})y_{k'}\\
{} & = & cy_{k}\sum_{l\geq k^{'}}X^{\mathbf{a},l}y_{l}-dy_{k^{'}}\sum_{l\geq k}X^{\mathbf{b},l}y_{l}-\sum_{t=1}^{\binom{n}{k}}\sum_{l\geq k}h_{t}X^{\mathbf{a}_{t},l}y_{l}y_{k^{'}} + T - T;
\end{eqnarray*}
where $T = c\sum_{l\geq k}[a_{1},\ldots, a_{k^{'}}\mid 1,\ldots,k-1,l,k+1,\ldots, k^{'}]y_{l}y_{k^{'}}$.  After a rearrangement of terms, we may write
\begin{eqnarray*}
r & = & \left(T-\sum_{t=1}^{\binom{n}{k}}\sum_{l\geq k}h_{t}X^{\mathbf{a}_{t},l}y_{l}y_{k^{'}} -dy_{k^{'}}\sum_{l\geq k}X^{\mathbf{b},l}y_{l}\right)\\
{} & {} & +\left(cy_{k}\sum_{l\geq k^{'}}X^{\mathbf{a},l}y_{l}\right) - T.\\[2mm]
\end{eqnarray*}
Let $T^{'}= c\sum_{l> k}[a_{1},\ldots, a_{k^{'}}\mid 1,\ldots,k-1,l,k+1,\ldots, k^{'}]y_{l}y_{k^{'}}$. Now we note,  $cX^{\textbf{a}}-dX^{\textbf{b}}- \sum_{t=1}^{\binom{n}{k}}h_{t}X^{\mathbf{a}_{t}}=0 $.
Hence $T-\sum_{t=1}^{\binom{n}{k}}\sum_{l\geq k}h_{t}X^{\mathbf{a}_{t},l}y_{l}y_{k^{'}} -
dy_{k^{'}}\sum_{l\geq k}X^{\mathbf{b},l}y_{l}$ becomes equal to 
$$T^{'}-\sum_{t=1}^{\binom{n}{k}}\sum_{l> k}h_{t}X^{\mathbf{a}_{t},l}y_{l}y_{k^{'}} -dy_{k^{'}}\sum_{l> k}X^{\mathbf{b},l}y_{l}.$$ 
\medskip

We also have 
$cy_{k}\sum_{l\geq k^{'}}X^{\mathbf{a},l}y_{l} - T = cy_{k}\sum_{l>
k^{'}}X^{\mathbf{a},l}y_{l} - T^{'}$, since the term for $l=k^{'}$ in $cy_{k}\sum_{l\geq k^{'}}X^{\mathbf{a},l}y_{l}$ gets cancelled with the term appearing in $T$ for $l=k$. Hence we write  
\begin{eqnarray*}
r & = & \left(T^{'}-\sum_{t=1}^{\binom{n}{k}}\sum_{l> k}h_{t}X^{\mathbf{a}_{t},l}y_{l}y_{k^{'}} -dy_{k^{'}}\sum_{l> k}X^{\mathbf{b},l}y_{l}\right)_{1}\\
{} & {} & +\left(cy_{k}\sum_{l> k^{'}}X^{\mathbf{a},l}y_{l}\right)_{2} - T^{'}\\[2mm]
{}&= & \left(\ \right)_{1} + \left(\ \right)_{2}- T^{'}.
\end{eqnarray*}

Clearly, the expression $(\ )_{1}$ belongs to $\widetilde{I}_{k+1}$, by Lemma \ref{laplace}. 
We note that no term of $(\ )_{1}$ contains $y_{k}$. So also for $T^{'}.$ Hence, the leading term of $r$ is the leading term of $(\ )_{2}$. 
By an application of similar argument as above we see that the expression $(\ )_{2}$, 
after division by elements of $\widetilde{S}_{k}$, further reduces to 
\begin{eqnarray*}
-\left(\sum _{l> k^{'}}\sum_{s\geq k^{'}}c[a_{1},\ldots, a_{k'}|1,\ldots, k-1 ,s ,k+1 ,\ldots ,k'-1, l]y_{l}y_{s}\right)\\
= \quad -\left(\sum _{l> k^{'}}\sum_{s> k^{'}}c[a_{1},\ldots, a_{k'}|1,\ldots, k-1 ,s ,k+1 ,\ldots ,k'-1, l]y_{l}y_{s}\right)\\
-\left(\sum _{l> k^{'}}c[a_{1},\ldots, a_{k'}|1,\ldots, k-1 ,k^{'} ,k+1 ,\ldots ,k'-1, l]y_{l}y_{k^{'}}\right).
\end{eqnarray*}
Moreover,
$$\sum _{l> k^{'}}c[a_{1},\ldots, a_{k'}|1,\ldots, k-1 ,k^{'} ,k+1 ,\ldots ,k'-1, l]y_{l}y_{k^{'}} + T' = 0$$ 
and
$$\sum _{l> k^{'}}\sum_{s>k^{'}}c[a_{1},\ldots, a_{k'}|1,\ldots, k-1 ,s ,k+1 ,\ldots ,k'-1, l]y_{l}y_{k^{'}}=0.$$ 
Therefore, after division by elements of $\widetilde{S}_{k}$, the expression 
$(\ )_{1}+(\ )_{2}-T'$ reduces to $(\ )_{1}$, which is in $\widetilde{I}_{k+1}$. \qed
\medskip 

\noindent\textbf{Proof of Theorem \ref{gbtheorem}.} 
We use induction on $n-k$ to prove that $G_{k}$ is a Gr\"{o}bner basis for the 
ideal $\widetilde{I}_{k}$. For $n-k=0$; the set $G_{k} = \widetilde{S}_{n}$ contains only one element 
and hence trivially forms a Gr\"{o}bner basis. We apply Buchberger's algorithm to 
prove our claim. Let $X^{\mathbf{a}},X^{\mathbf{b}} \in G_{k}$. The following 
cases may arise: 
\begin{itemize}
\item $X^{\mathbf{a}},X^{\mathbf{b}}\in S_{k}$, for $\mathbf{a}, \mathbf{b}\in C_{k}$; 
\item $X^{\mathbf{a}}\in S_{k'}$ and $X^{\mathbf{b}}\in S_{k}$ where $k'>k$; $\mathbf{a}\in C_{k'}$ and $\mathbf{b}\in C_{k}$.
\end{itemize}
We have proved in Lemmas \ref{spolylemma} and \ref{spolydifferent} that upon 
division by $\widetilde{S}_{k}$, the $S$-polynomial 
$S(\widetilde{X}^{\mathbf{a}},\widetilde{X}^{\mathbf{b}})\longrightarrow r$ 
for some $r\in\widetilde{I}_{k+1}$, in both the cases. 
By induction hypothesis, $G_{k+1}$ is a Gr\"{o}bner basis for $\widetilde{I}_{k+1}$. 
Hence $r$ reduces to $0$ modulo $G_{k+1}$ and hence modulo $G_{k}$, 
since $G_{k+1}\subset G_{k}$ .
\medskip

We now show that $G_{k}$ is a reduced Gr\"{o}bner basis for $\widetilde{I}_{k}$. 
Let $X^{\mathbf{a}}\in S_{k'}$ and $X^{\mathbf{b}}\in S_{k}$ where $k'\geq k$; 
$\mathbf{a}\in C_{k'}$ and $\mathbf{b}\in C_{k}$. Then, 
$\widetilde{X}^{\mathbf{a}} = \sum_{i\geq k'}X^{\mathbf{a},i}y_{i}$ and 
$\widetilde{X}^{\mathbf{b}} = \sum_{i\geq k}X^{\mathbf{b},i}y_{i}$. 
If $k' > k$, then $y_{k'}|\LT(\widetilde{X}^{\mathbf{a}})$ but does not divide 
$\LT(\widetilde{X}^{\mathbf{b}})$. Hence, $\LT(\widetilde{X}^{\mathbf{a}})$ does not 
divide $\LT(\widetilde{X}^{\mathbf{b}})$. If $k'=k$, then 
$\LT(\widetilde{X}^{\mathbf{a}})=x_{(a_{1},1)}\cdots x_{(a_{k},k)}y_{k}$ and $\LT(\widetilde{X}^{\mathbf{b}})=x_{(b_{1},1)}\cdots x_{(b_{k},k)}y_{k}$. Therefore, 
$\widetilde{X}^{\mathbf{a}}|\widetilde{X}^{\mathbf{b}}$ implies that 
$\textbf{a}=\textbf{b}$. This proves that the Gr\"{o}bner basis is reduced. \qed

\section{Gr\"{o}bner basis for $\mathcal{J}$}

\begin{theorem}\label{gbtheoremJ}
Let us consider the lexicographic monomial order induced by 
$y_{1}>y_{2}>\cdots >y_{n}>x_{11}>x_{12}>\cdots > x_{(n+1),(n-1)}>x_{(n+1),n}$ 
on $\widehat{R} = K[x_{ij}, y_{j}\mid 1\leq i \leq n+1, \, 1\leq j \leq n]$. 
The set $G_{k}$ is a reduced Gr\"{o}bner Basis for the ideal 
$\widetilde{I}_{k}$. In 
particular, $\mathcal{G} = G_{1}$ is a reduced Gr\"{o}bner Basis for the ideal 
$\widetilde{I}_{1} = \mathcal{J}$. 
\end{theorem}
\proof The scheme of the proof is the same as that for $\mathcal{I}$, with suitable 
changes made for $\widehat{X}$ in the Lemmas. We only reiterate the last part of 
the proof where we carry out induction on $n-k$. For $n-k=0$, the set 
$G_{k} = \widetilde{S}_{n} = \{\Delta_{1}y_{n}, \ldots , \Delta_{n+1}y_{n}\}$, where 
$\Delta_{i} = \det(\widehat{X}_{i})$. We first note that $\LT(\Delta_{i})$ and 
$\LT(\Delta_{j})$ are coprime. Therefore, 
\begin{eqnarray*}
S(\Delta_{i}y_{n}, \Delta_{j}y_{n}) & = & \LT(\Delta_{j})\cdot(\Delta_{i}y_{n}) - \LT(\Delta_{i})\cdot(\Delta_{j}y_{n})\\
{} & = & \LT(\Delta_{j})(\LT(\Delta_{i})y_{n} + y_{n}p_{i}) - 
\LT(\Delta_{i})(\LT(\Delta_{j})y_{n} - y_{n}p_{j})\\
{} & = & (\LT(\Delta_{j})y_{n})p_{i} - (\LT(\Delta_{i})y_{n})p_{j}\\
{} & = & (\Delta_{j}y_{n} - p_{j}y_{n})p_{i} - (\Delta_{i}y_{n} - p_{i}y_{n})p_{j}\\
{} & = & \Delta_{j}y_{n}p_{i} - \Delta_{i}y_{n}p_{j}\longrightarrow _{G_{n}} 0.
\end{eqnarray*}
The rest of the proof is essentially the same as that for Theorem \ref{gbtheorem}.\qed

\section{Betti Numbers of $\mathcal{I}$ and $\mathcal{J}$}
\begin{theorem}\label{regseqI}
Suppose that $X = (x_{ij})_{n\times n}$ is either a generic or a generic symmetric 
$n\times n$ matrix and $Y$ a generic $n\times 1$ matrix given by $Y=(y_{j})_{n\times 1}$. If $X$ is generic, we write $g_{i}=\sum_{j=1}^{n} x_{ij}y_{j}$ and $\mathcal{I} = I_{1}(XY) = \langle g_{1},g_{2},\cdots ,g_{n}\rangle$. If $X$ is 
generic symmetric, we write $g_{1}=\sum_{j=1}^{n}x_{1j}y_{j}$, $g_{n}=(\sum_{1\leq k \leq n}x_{kn}y_{k})$ and $g_{i}=(\sum_{1\leq k<i}x_{ki}y_{k})+(\sum_{i\leq k\leq n}x_{ik}y_{k})$ for $1<i<n$ and $\mathcal{I} = I_{1}(XY) = \langle g_{1},\cdots ,g_{n}\rangle$. The generators $g_{1}, \ldots , g_{n}$ of $\mathcal{I} = I_{1}(XY)$ in either 
case form a regular sequence in the polynomial $K$-algebra $R = K[x_{ij}, \, y_{j} \mid 1\leq i,j\leq n]$. Moreover, $\{g_{1}, \ldots , g_{n}\}$ form a Gr\"obner basis for 
$\mathcal{I}$ in either case with respect to the lexicographic monomial order which satisfies (1) and (2) given below:
\begin{enumerate}
\item $x_{11}> x_{22}> \cdots >x_{nn}$;
\item $x_{ij}, y_{j} < x_{nn}$ for every $1 \leq i \neq j \leq n$.
\end{enumerate}
\end{theorem} 
\medskip

\proof The monomial order chosen is lexicographic order induced by the ordering among the variables given by (1) and (2). It is clear from the expressions of $g_{i}$ that their leading terms are pairwise coprime. Therefore, the proof follows from Lemma \ref{disjoint}. \qed
\medskip

\begin{corollary}\label{bettiI}
$\mathcal{I}$ is minimally resolved by the Koszul complex $\mathbb{G}$ and the $i$-th Betti number 
of $\mathcal{I}$ is $\binom{n}{i}$.
\end{corollary}  
\medskip

\begin{theorem}\label{bettiJ}
Suppose that $\widehat{X} = (x_{ij})_{(n+1)\times n}$ is a generic $(n+1)\times n$ 
matrix and $Y$ a generic $n\times 1$ matrix given by $Y=(y_{j})_{n\times 1}$. 
Let $g_{i}=\sum_{j=1}^{n+1} x_{ij}y_{j}$ and $\mathcal{J} = I_{1}(\widehat{X}Y) = 
\langle g_{1},\cdots ,g_{n+1}\rangle$. The total Betti numbers of the ideal 
$\mathcal{J}$ are $\beta_{0}=1, \beta_{1}=n+1$, $\beta_{n+1}= n$, 
$\beta_{k+1}=\binom{n}{k}+\binom{n}{k-1}+\binom{n}{k+1}$ for $1\leq k < n$.
\end{theorem} 
\medskip

We first discuss the scheme of the proof below. 
We will use the following observations to compute the 
total Betti numbers of $\mathcal{J}$.
\begin{itemize}
\item[Step 1.] The minimal graded free resolution of 
$\mathcal{I}=\langle g_{1},\cdots,g_{n}\rangle $ is given by the Koszul Resolution.

\item[Step 2.] We prove that $\langle g_{1},\cdots,g_{n}: g_{n+1}\rangle = \langle g_{1},\cdots,g_{n},\Delta\rangle$; where $\Delta= \det(X)$. This proof requires the fact that $\langle g_{1},\cdots,g_{n},\Delta\rangle$ is a prime ideal, which 
has been proved in Theorem 5.4 in \cite{sstprime}. 

\item[Step 3.] We prove that $\langle g_{1},\cdots g_{n}:\Delta\rangle =\langle y_{1},y_{2},\cdots, y_{n}\rangle$.

\item[Step 4.] We construct a graded free resolution of 
$\langle g_{1},\cdots,g_{n},\Delta\rangle $ 
using mapping cone between resolutions of $\langle g_{1},\cdots,g_{n}\rangle$ 
and $\langle y_{1},\cdots,y_{n}\rangle $. We extract a minimal free resolution 
from this resolution. 

\item[Step 5.] Finally, we construct a graded free resolution of 
$\langle g_{1},\cdots,g_{n},g_{n+1}\rangle$ using mapping cone between free 
resolutions of $\langle g_{1},\cdots,g_{n},\Delta\rangle$ and 
$\langle g_{1},\cdots,g_{n}\rangle$. We extract a minimal free resolution 
from this resolution.
\end{itemize}
\medskip

\begin{remark}
We need detailed information about the 
ideal $\langle g_{1},\cdots,g_{n},\Delta\rangle$, where $\Delta = \det(X)$. We need the fact that this ideal is a 
prime ideal, which has been proved in Theorem 5.4 in \cite{sstprime}. We also need a 
minimal free resolution for this ideal, which has been proved below in Lemma \ref{resgennorth}. 
We came to know much later that $\langle g_{1},\cdots,g_{n},\Delta\rangle$ was defined in \cite{nor}. 
It is known as the generic Northcott ideal and a minimal free resolution can be found in \cite{nor}. However, 
we give a different proof here using our Gr\"{o}bner basis computation, which also shows the linking 
of nested complete intersection ideals. Moreover, Northcott's resolution can perhaps be used to prove that 
$\langle g_{1},\cdots,g_{n},\Delta\rangle$ is a prime ideal, although our proof in \cite{sstprime} is 
absolutely different and uses the result in \cite{fer}.
\end{remark}
\medskip

\begin{lemma}\label{cofactor}
$\Delta y_{i}=\sum_{j=1}^{n}A_{ji}g_{j}$, where $A_{ji}$ is the cofactor of 
$x_{ji}$ in $X$.
\end{lemma}

\proof We have 
$$\Delta y_{i} = \sum_{j=1}^{n}A_{ji}x_{ji}y_{i} 
 = \sum_{j=1}^{n}A_{ji}\left(\sum_{k=1}^{n}x_{jk}y_{k}\right) 
 - \sum_{j=1}^{n}A_{ji}\left(\sum_{k\neq i}x_{jk}y_{k}\right) 
 = \sum_{j=1}^{n}A_{ji}g_{j},$$ 
 since $\sum_{j=1}^{n}A_{ji}\left(\sum_{k\neq i}x_{jk}y_{k}\right) 
 = \sum_{k\neq i}\left(\sum_{j=1}^{n}A_{ji}x_{jk}\right)y_{k} = 0$.\qed
\medskip

\begin{lemma}\label{cofactorcor} $\langle g_{1},\cdots,g_{n},\Delta\rangle\subseteq \langle g_{1},\cdots,g_{n}: g_{n+1}\rangle$.
\end{lemma}

\proof We have 
$g_{i}\in \langle g_{1},\cdots,g_{n}: g_{n+1}\rangle$, for every 
$1\leq i\leq n$. Moreover, $y_{i}\Delta\in \langle g_{1},\cdots,g_{n}\rangle$, 
by Lemma \ref{cofactor}. Hence, 
$g_{n+1}\Delta\in \langle g_{1},\cdots,g_{n}\rangle$. \qed
\medskip

\begin{lemma}\label{colon1}
$\langle g_{1},\cdots,g_{n}: g_{n+1}\rangle=\langle g_{1},\cdots,g_{n},\Delta\rangle $
\end{lemma}

\proof We have proved that $\langle g_{1},\cdots,g_{n},\Delta\rangle\subseteq \langle g_{1},\cdots,g_{n}: g_{n+1}\rangle$ in Lemma \ref{cofactorcor}. We now prove that 
$\langle g_{1},\cdots,g_{n}:g_{n+1}\rangle\subseteq \langle g_{1},\cdots,g_{n},\Delta\rangle$. Let $z\in \langle g_{1},\cdots,g_{n}:g_{n+1}\rangle$. Then 
$zg_{n+1}\in \langle g_{1},\cdots,g_{n}\rangle\subset\langle g_{1},\cdots,g_{n},\Delta\rangle$. It is easy to see that $g_{n+1}\notin \langle g_{1},\cdots,g_{n},\Delta\rangle$. Therefore, $z\in \langle g_{1},\cdots,g_{n},\Delta\rangle$, since 
$\langle g_{1},\cdots,g_{n},\Delta\rangle$ is a prime ideal by Theorem 5.4 in \cite{sstprime}. \qed
\medskip

\begin{lemma}\label{colon2}
$\langle g_{1},\cdots,g_{n}:\Delta\rangle=\langle y_{1},\cdots,y_{n}\rangle$
\end{lemma}

\proof We have $y_{i}\Delta\in \langle g_{1},\cdots,g_{n}\rangle$ by Lemma \ref{cofactor}; which implies that $\langle y_{1},\cdots,y_{n}\rangle\subset\langle g_{1},\cdots,g_{n}:\Delta\rangle$. Let 
$z\in \langle g_{1},\cdots,g_{n}:\Delta\rangle$. Then 
$z\Delta\in\langle g_{1},\cdots,g_{n}\rangle\subseteq \langle y_{1},\cdots,y_{n}\rangle$. Therefore, $z\in \langle y_{1},\cdots,y_{n}\rangle$, since 
$\Delta\notin \langle y_{1},\cdots,y_{n}\rangle$ and 
$\langle y_{1},\cdots,y_{n}\rangle$ is a prime ideal.\qed
\medskip

\noindent\textbf{Mapping Cones.} \, The resolution for $\langle y_{1},\cdots,y_{n}\rangle$ is given by the Koszul complex 
$\mathbb{F}_{\centerdot}$\, . We now give a resolution of 
$\langle g_{1},\cdots,g_{n},\Delta\rangle$ by the mapping cone technique. 
We know that $\langle g_{1},
\cdots,g_{n}:\Delta\rangle=\langle y_{1},\cdots,y_{n}\rangle$, by 
Lemma \ref{colon2}. We first construct a connecting homomorphism 
$\phi_{\centerdot}:\mathbb{F}_{\centerdot}\longrightarrow \mathbb{G}_{\centerdot}$\, . 
Let $\phi_{0}$ denote the multiplication by $\Delta$. In order to make 
the map $\phi_{0}$ a degree zero map, we set the grading as 
$\mathbb{F}_{0}\cong (R(-n))^{1}$ and $\mathbb{G}_{0}= (R(0))^{1}$. 
Since $\mathbb{F}_{\centerdot}$ and $\mathbb{G}_{\centerdot}$ are both 
Koszul resolutions, we set the grading as $\mathbb{G}_{i}\cong (R(-2i))^{\binom{n}{i}}$ 
and $\mathbb{F}_{i}\cong (R(-n-i))^{\binom{n}{i}}$. Now we see that, 
$i\neq n$ implies that $-2i\neq -n -i$. Hence the image of $\phi_{i}$ 
for $i\neq n$ is contained in the maximal ideal. We have 
$\mathbb{F}_{i}=\mathbb{G}_{i}$, only for $i=n$. If we can show 
that the map $\phi_{n}$ is not the zero map, then this will be the 
only free part of the resolution which we can cancel out for obtaining the 
minimal resolution.
\medskip

\begin{lemma}
The map $\phi_{n}$ is not the zero map.
\end{lemma}

\proof We refer to  \cite{hip}. If $\phi_{n}$ is the zero map, then $\phi_{0}(R)\subseteq \delta_{1}(\mathbb{G}_{1})$, where $\delta_{.}$ denotes the differential of $\mathbb{G}$. The image of $\delta_{1}$ is the ideal $\langle g_{1},\cdots,g_{n}\rangle$, which does not 
contain $\phi_{0}(1)=\Delta$. The map $\phi_{n}$ is not the zero map.\qed
\medskip

Therefore, the above discussion proves the following Lemma.
\medskip

\begin{lemma}\label{resgennorth}
Hence a minimal graded free resolution of $\langle g_{1},\cdots,g_{n},\Delta\rangle$ is given by $\mathbb{M_{\centerdot}}$, such that $\mathbb{M}_{i}\cong (R(-n-i+1))^{\binom{n}{i-1}}\oplus (R(-2i))^{\binom{n}{i}}$ for $0<i<n$, $\mathbb{M}_{0}\cong R(0)$ and 
$\mathbb{M}_{n}\cong (R(-2n))^{n}$.
\end{lemma}
\medskip

\noindent\textbf{(Proof of Theorem \ref{bettiJ}.)} We now find the Betti numbers for the ideal $\langle g_{1},\cdots, g_{n+1}\rangle$ by constructing the mapping cone between the resolutions $\mathbb{M_{\centerdot}}$ and the resolution $\mathbb{G_{\centerdot}}$ of 
$\langle g_{1},\cdots, g_{n}\rangle$. The connecting map $\psi_{0}$ is multiplication 
by $g_{n+1}$. Hence to make it degree zero we set, $\mathbb{G}_{0}= (R(2))^{1}$ and 
$\mathbb{G}_{i}\cong (R(2-2i))^{\binom{n}{i}}$ for $i>0$. Here we note that 
$2-2i\neq -2i$ and $-n-i+1\neq 2-2i$ for $1\leq i\leq n$. Hence, for each $1\leq i\leq n$, the image of $\psi_{i}$ is contained in the maximal ideal. This shows that the 
resolution obtained by the mapping cone between $\mathbb{M_{\centerdot}}$ and 
$\mathbb{G_{\centerdot}}$ is minimal. Hence the total Betti numbers of $\mathcal{J}$ 
are:
\medskip

\noindent $\beta_{0}=1, \beta_{1}=n+1$; \\
$\beta_{n+1}= n$; \\
$\beta_{k+1}=\binom{n}{k}+\binom{n}{k-1}+\binom{n}{k+1}$ for $1\leq k < n$. \qed

\medskip

\begin{corollary}
The ring $R/\mathcal{I}$ is Cohen-Macaulay and the ring $\hat{R}/\mathcal{J}$ is not 
Cohen-Macaulay.
\end{corollary}
\proof
The polynomial ring $R$ is Cohen-Macaulay and $g_{1},\ldots,g_{n}$ is a regular sequence therefore the ring $R/\mathcal{I}$ is Cohen-Macaulay.
\medskip

We have seen that $\textrm{projdim}_{\widehat{R}}\widehat{R}/\mathcal{J}=n+1$. 
Therefore, by the Auslander-Bauchsbaum formula 
$\textrm{depth}_{\widehat{R}}\widehat{R}/\mathcal{J}=n(n+1)+n-(n+1)=n^{2}+n-1$. 
We have proved in Lemma 5.5 in \cite{sstprime} that 
$\langle y_{1},\ldots,y_{n}\rangle$ is a minimal prime over $\mathcal{J}$. 
Therefore, $\textrm{dim}\widehat{R}/\mathcal{J}\geq \textrm{dim}\widehat{R}/\langle y_{1},\ldots,y_{n}\rangle=n^{2}+n $; hence the ring $\widehat{R}/\mathcal{J}$ 
is not Cohen-Macaulay.\qed

\section{$I_{1}(XY)$, where X is $m\times mn$ generic matrix and Y is 
$mn\times n$ generic matrix}

Finally, we consider the case when $X=(x_{ij})_{m\times mn}$ is a 
generic matrix of size $m\times mn$ and $Y=(y_{ij})_{mn\times n}$ 
is generic matrix of size $mn\times n$. We define 
$\mathfrak{I} = I_{1}(XY)$. Let $g_{ij} = \sum_{t=1}^{mn}x_{it}y_{tj}$, 
with $1\leq i\leq m, \, 1\leq i\leq n$. Then, 
$\mathfrak{I} = \langle \{g_{ij} \mid 1\leq i\leq m, \, 1\leq i\leq n\}\rangle$. In this section we construct a Gr\"{o}bner basis for the ideal 
$\mathfrak{I}$ with respect to a suitable monomial order and use 
that to show that the generators $g_{ij}$, 
with $1\leq i\leq m$, 
$1\leq i\leq n$ form a regular sequence. We first set a few notations 
before we prove the main results. 
\begin{itemize}
\item $X=\begin{pmatrix}A_{1}&\cdots & A_{n}
\end{pmatrix}$, where $A_{s}=\begin{pmatrix}
x_{1(m(s-1)+1)}&\cdots & x_{1(ms)}\\
\vdots &  \vdots&  \vdots\\
x_{m(m(s-1)+1)}&\cdots & x_{m(ms)}\\
\end{pmatrix} $ is the $m\times m$ matrix for every $1\leq s\leq n$.
\medskip

\item $[X]_{s}=\begin{pmatrix}
A_{s}& A_{1}&\cdots &\widehat{A_{s}}&\cdots & A_{n}
\end{pmatrix}$, for every $1\leq s\leq n$.
\medskip

\item $[Y]_{s}=\begin{pmatrix}
y_{(m(s-1)+1)s}\\
\vdots\\
y_{(ms)s}\\
y_{1s}\\
\vdots\\
y_{(mn)s}
\end{pmatrix}$, for every $1\leq s\leq n$.
\end{itemize}
\medskip

We will use Theorem \ref{gbtheorem} for constructing a Gr\"{o}bner basis for the 
ideal $\mathfrak{I}$. A very important reason behind considering this 
class of ideals is that we get some nice examples of transversal intersection 
of ideals. Two results that would be useful for our purpose are the following:
\medskip
 
\begin{lemma} \label{trans}
Let $>$ be a monomial ordering on $R$. Let $I$ and $J$ be ideals in $R$, such that 
$m(I)$ and $m(J)$ denote unique minimal generating sets for their leading ideals 
$\LT(I)$ and $\LT(J)$ respectively. Then, $I\cap J = IJ$ if the set of variables 
occurring in the set $m(I)$ is disjointed from the the set of variables 
occurring in the set $m(J)$.
\end{lemma}

\proof See Lemma 3.6 in \cite{sstsum}.\qed
\medskip

\begin{lemma}\label{tensorprod}
Let $I$ and $J$ be graded ideals in a graded ring $R$, such that $I\cap J=I\cdot J$. Suppose that $\mathbb{F}_{\centerdot}$ and $\mathbb{G}_{\centerdot}$ are minimal free resolutions of $I$ and $J$ respectively. Then $\mathbb{F}_{\centerdot}\otimes \mathbb{G}_{\centerdot}$ is a minimal free resolution for the graded ideal $I+J$.
\end{lemma}

\proof See Lemma 3.7 in \cite{sstsum}.\qed
\medskip 

\begin{theorem}\label{gmain}
Let us choose the lexicographic monomial order on $R$ induced by $ y_{11}>y_{21}>\cdots >y_{(mn)1}> y_{(m+1)2}>y_{(m+2)2}>\cdots >y_{(2m)2}>y_{12}>\cdots y_{(mn)2}>\cdots > y_{(m(n-1)+1)n}>y_{(m(n-1)+2)n}>\cdots >y_{((mn)n}>y_{1n}>\cdots y_{(m(n-1))n}>x_{11}> x_{12}> \cdots >x_{m(mn)}$. Let $\mathcal{G}_{s}$ be the reduced Gr\"{o}bner Basis of the ideal  $I_{1}([X]_{s}[Y]_{s})$ for $1\leq s\leq n$,  obtained by Theorem \ref{gbtheorem}. Then $\mathfrak{G}_{t}=\cup_{s=1}^{t}\mathcal{G}_{s}$ is a reduced Gr\"{o}bner Basis for the ideal  $P_{t}=\sum_{s=1}^{t}I_{1}([X]_{s}[Y]_{s})$ for $1\leq t\leq n$. In particular, $\mathfrak{G}_{n}$ is a reduced Gr\"{o}bner Basis 
for the ideal $P_{n}=\mathfrak{I}=I_{1}(XY)$.
\end{theorem}

\proof We have $P_{t}=\sum_{s=1}^{t}I_{1}([X]_{s}[Y]_{s})$, and we observe 
that if $p\in \mathcal{G}_{s} $ and $q\in \mathcal{G}_{t}$ for $1\leq s<t\leq n$, 
then $\gcd(\LT(p),\LT(q))=1$. Therefore the $S$-polynomial of $p,q$ reduces to 
zero after applying division upon $\mathfrak{G}_{t}$.\qed

\begin{theorem}
Let us denote $P_{t}=\sum_{s=1}^{t}I_{1}([X]_{s}[Y]_{s})$, for $1\leq t\leq n-1$. Then $P_{t}\cap I_{1}([X]_{t+1}[Y]_{t+1})= P_{t}\cdot I_{1}([X]_{t+1}[Y]_{t+1})$. Hence the elements $g_{ij} = \sum_{t=1}^{mn}x_{it}y_{tj}$, $1\leq i\leq m$, 
\, $1\leq i\leq n$ form a regular sequence and the 
Koszul complex resolves $R/\mathfrak{I}$ as an $R$-module minimally.
\end{theorem}

\proof If $p\in \mathcal{G}_{s} $ and $q\in \mathcal{G}_{t}$, for $1\leq s<t\leq n$. 
Then $\gcd(\LT(p),\LT(q))=1$, therefore by theorem \ref{gmain} and lemma \ref{trans}, we have $P_{t}\cap I_{1}([X]_{t+1}[Y]_{t+1})= P_{t}\cdot I_{1}([X]_{t+1}[Y]_{t+1})$.
\medskip

By Theorem \ref{regseqI} the generators of the ideal $P_{1}$ form a regular sequence and also the generators of the ideal $I_{1}([X]_{s}[Y]_{s})$ form a regular sequence for each $1\leq s\leq n$. Hence the Koszul complex resolve 
$R/P_{1}$ and $R/I_{1}([X]_{s}[Y]_{s})$ minimally. Now 
$P_{t}\cap I_{1}([X]_{t+1}[Y]_{t+1})= P_{t}\cdot I_{1}([X]_{t+1}[Y]_{t+1})$. Hence, by application of lemma \ref{trans} we can conclude that the Koszul 
complex resolves $R/\mathfrak{I}$ minimaly. \qed

\section*{Acknowledgements}
The second author is the corresponding author who has been supported by 
the research project EMR/2015/000776, sponsored by the SERB, Government of 
India. The third author thanks SERB for the post-doctoral fellowship under 
the said project. The third author thanks CSIR for the Senior Research 
Fellowship for Ph.D. The authors thank the anonymous referees for their 
valuable comments and for drawing their attention to the references 
\cite{concini} and \cite{tchernev}, extremely pertinent to this work.

\bibliographystyle{amsalpha}

\end{document}